\documentclass{article}
\usepackage[top=2cm, bottom=2cm, left=3cm, right=3cm]{geometry}
\usepackage{amssymb}
\usepackage{amsmath}
\usepackage{latexsym}
\usepackage{mathrsfs}
\usepackage{graphicx}
\usepackage[all]{xy}
\usepackage{hyperref}
\newtheorem{Th}{Theorem}
\newtheorem{Lemma}[Th]{Lemma}

\title{Puzzle Model for Bumpless Pipe Dream}
\author{XIONG Rui}
\begin{document}

\maketitle

\begin{abstract}
  In this note, a new puzzle is introduced where the pipe dream and bumpless pipe dream can be played simultaneously.
  Using these, a combinatorial proof of the (ordinary) Schubert polynomials in terms of bumpless pipe dream is given.
  The main tool is the Yang--Baxter equation.
\end{abstract}
I would thank for Paul Zinn-Justin, 
Evgeny Smirnov and Maksim Karev
for their encouragement.

\section{Introduction}

\paragraph{Schubert Polynomials. }%
Consider the polynomial ring $R$ over $\mathbb{Q}$ in indeterminants $x_1,\ldots,x_n,\ldots$. 
We define the \emph{Demazure operator} over it, for $i=1,2,\ldots$
$$\partial_i f:=\frac{f(x)-f(s_ix)}{x_i-x_{i+1}},$$
where $s_i=(i,i+1)\in \mathfrak{S}_\infty$ a simple reflection of infinite symmetric group.
In \cite{lascoux1982structure}, Lascoux and Sch{\"u}tzenberger defined the \emph{Schubert polynomials}.
by the relation
$$\mathfrak{S}_{w_0^n}(x)=x_1^{n-1}\ldots x_{n-1},\qquad \partial_i\mathfrak{S}_{w}=\begin{cases}
\mathfrak{S}_{ws_i},& \ell(ws_i)=\ell(w)-1, \\
0, &\text{otherwise},
\end{cases}$$
where $w_0^n$ is the longest word in $\mathfrak{S}_n$, and $\ell$ is the length function. 

\paragraph{Pipe Dreams. }%
\def\PD{\mathsf{PD}}
\def\wt{\operatorname{\mathtt{wt}}}
In \cite{bergeron1993rc}, a combinatorial description of Schubert polynomials was discovered, which is now known as \emph{pipe dream}.
It is a tiling of $1/4$ plane by two kinds of pieces with no pair of pipes crossing more than twice.
Here is an example.
$$\includegraphics{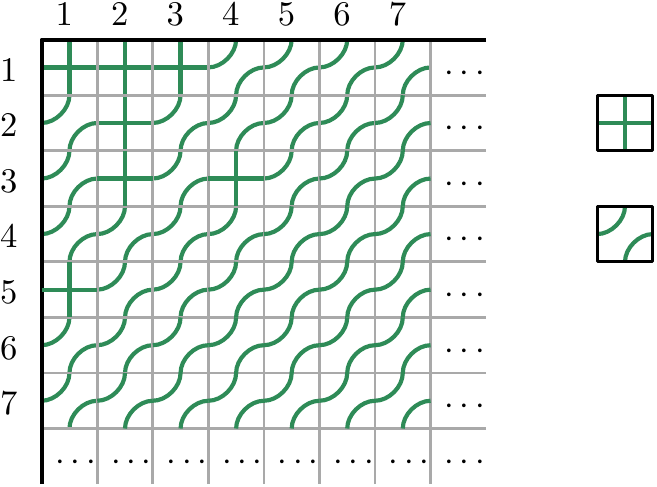}$$
For any pipe dream $\alpha$, it determines a permutation $w(\alpha)$,
say, $w(i)=j$ if the leftmost $i$-th pipe is connected to the upper $j$-th pipe.
For any permutation $w$, we denote $\PD(w)$ all pipe dreams $\alpha$ with $w(\alpha)=w$.
We define its weight $\wt(\alpha)$ to be the product of all $x_i$ if there is a cross in $(i,j)$-position.
It is known that

\begin{Th}[\cite{bergeron1993rc}]\label{ThPD}For any $w\in \mathfrak{S}_\infty$, 
$\mathfrak{S}_w(x)=\sum_{\alpha\in \PD(w)} \wt(\alpha)$. 
\end{Th}

This combinatorial picture is essentially an expansion of the generating function of Schubert polynomials
discovered in \cite{fomin1994schubert}, see also \cite{fomin1996yang}.

\paragraph{Bumpless Pipe Dreams. }%
\def\BPD{\mathsf{BPD}}
In \cite{lam2018stable}, a new pipe dream was found, called the \emph{bumpless pipe dream}.
It is a tilting of $n\times n$ square with $6$-patterns with no pair of pipes crossing more than twice.
Here is an example.
$$\includegraphics{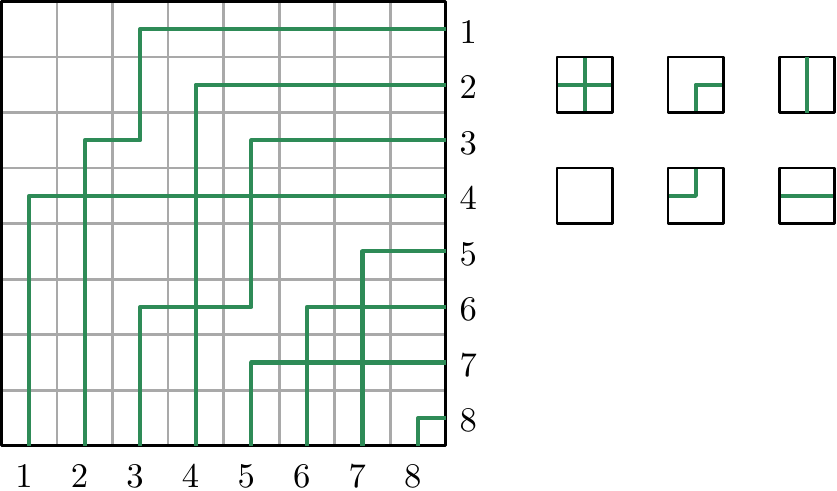}$$
For any bumpless pipe dream $\alpha$, it determines a permutation $w(\alpha)$,
say, $w(i)=j$ if the rightmost $i$-th pipe is connected to the lower $j$-th pipe.
For any permutation $w$, we denote $\BPD(w)$ all bumpless pipe dreams $\alpha$ with $w(\alpha)=w$.
We define its weight $wt(\alpha)$ to be the product of all $x_i$ if there is an blank in $(i,j)$-position.
For $w\in \mathfrak{S}_\infty$, it was proved in \cite{lam2018stable} in a long geometric way that

\begin{Th}[\cite{lam2018stable}]\label{ThBPD}For any $w\in \mathfrak{S}_\infty$, $\mathfrak{S}_w(x)=\sum_{\alpha\in \BPD(w)} \wt(\alpha)$. 
\end{Th}

Actually, in \cite{lam2018stable}, Theorem \ref{ThBPD} is proved for double Schubert polynomials (see below).

\paragraph{Double Schubert Polynomials. }%
We can define \emph{double Schubert polynomial} or \emph{equivariant Schubert polynomial}
by exchanging $\mathfrak{S}_{w_0^n}(x)$ by $\mathfrak{S}_{w_0^n}(x,y)=\prod_{i+j\leq n}(x_i-y_{j})$
in the definition of the Schubert polynomials. 

Meanwhile, in pipe dream, 
if we exchange $\wt(\alpha)$ by the product of all $(x_i-y_j)$ if there is a cross in $(i,j)$-position, 
then we still have $\mathfrak{S}_w(x,y)=\sum_{\alpha\in \PD(w)} \wt(\alpha)$ for all $w\in \mathfrak{S}_{\infty}$, 
see \cite{knutson2019schubert}.
In bumpless pipe dream, 
if we exchange $\wt(\alpha)$ by the product of all $(x_i-y_j)$ if there is an blank in $(i,j)$-position, 
then it is proved in \cite{lam2018stable} that $\mathfrak{S}_w(x,y)=\sum_{\alpha\in \BPD(w)} \wt(\alpha)$. 

It is also known by examples that summands of above two expressions are far from being bijective,
so a combinatorial proof is desired.

\paragraph{Notes. }%
The main purpose of this note is to show Theorem \ref{ThBPD} in a pure combinatorial way,
assuming Theorem \ref{ThPD}.
The proof uses a new puzzle where the pipe dream and bumpless pipe dream can be played simultaneously.
The main tool is the Yang--Baxter equation.

The author did not figure out how to generalize this proof to the double case.

\section{Puzzles}

Now, consider the following puzzles.
$$\includegraphics{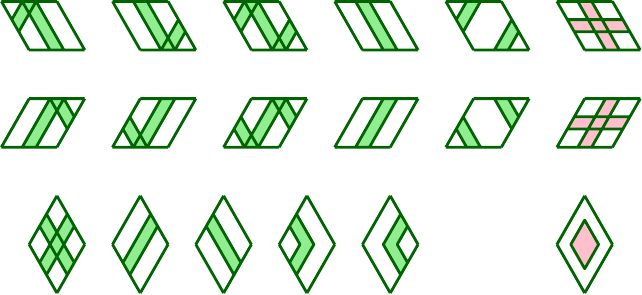}$$
Here are some terminologys of these puzzles
\begin{itemize}
  \item The inner pattern is referred to as \emph{pipes}.
  \item The puzzle in pink is said to be \emph{valued}.
\end{itemize}
Note that, the last puzzle has an empty pattern, but to distinguish it from an unknown puzzle, we draw a diamond inside.
To play with them, we define the following concepts
\begin{itemize}
  \item A \emph{chess board} is a certain tiling of three directions of parallelograms in some shape.
  \item A chessboard has a \emph{valuation} if each parallelogram is valued by some element in some commutative ring.
  \item A \emph{pipe tiling} of a chessboard is a pose of puzzles in a parallelogram with all pipes connected, and no pair of pipes crossing twice.
  Usually, it is denoted by letters around the shape.
  \item A \emph{rule} for a shape if the start and the end of pipes are given on the boundary.
  \item A \emph{solution} to a \emph{rule} is a pipe tiling satisfying the rule.
  \item The \emph{value} of a tiling is the product of valuation of diamond posed by a valued puzzle.
  \item The \emph{value} of a rule is the sum of values of all solutions.
\end{itemize}
For example,
$$\includegraphics{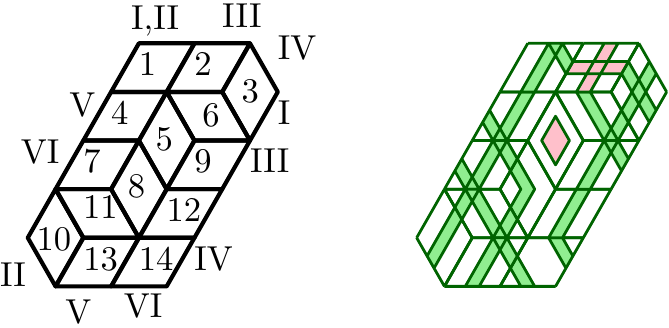}$$
the left hand is a valued, ruled chessboard, and the right-hand side the a solution with value $2+5$.

Here are some observation of our puzzles
\begin{itemize}
  \item We can put an orientation on pipes, with all pipes never going down,
  and horizontally left to right depending on the direction of the parallelogram.
  $$\includegraphics{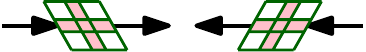}$$
  Note that any pipe tiling is orientation compactible.
  In particular, there will be no circle in any pipe tilting.
  \item For any horizontal side, the number of pipes can only be one or two;
  for the tilted side, the number of pipes can only be one or zero.
  Furthermore, if we denote the number of pipes around a puzzle as the following,
  $$\includegraphics{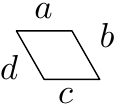}$$
  then $a+b=c+d$.
  In particular, $a+b\leq 2$ if and only if $c+d\leq 2$.
\end{itemize}

The relation of puzzles and pipe dreams are described in the following theorem.

\begin{Th}\label{puzzlePD}For a permutation $w\in \mathfrak{S}_n$,
then the value of the following chess board is the double Schubert polynomial $\mathfrak{S}_w(x,y)$ computed by pipe dream.
$$\includegraphics{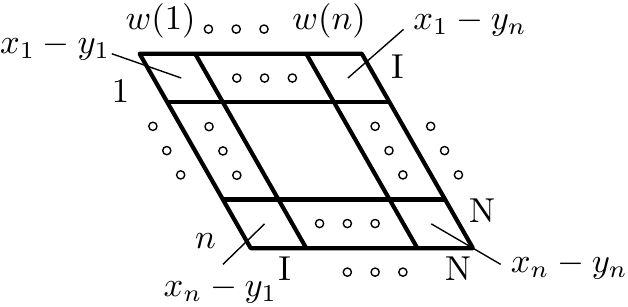}$$
\end{Th}
\noindent\textsc{Proof. }
Consider each row, there is one pipe going in and one going out in total,
so any horizontal side cannot have two pipes.
By intuition, for each parallelogram, there must be two pipes, so only two puzzles can be used.
It is obviously the same as the pipe dream. \qquad Q.E.D.
\bigbreak

This is a special case of Lemma \ref{YoungDiagram2} below.

\begin{Th}\label{puzzleBPD}For a permutation $w\in \mathfrak{S}_n$,
then the value of the following chess board is the double Schubert polynomial $S_w(x,y)$ computed by bumpless pipe dream.
$$\includegraphics{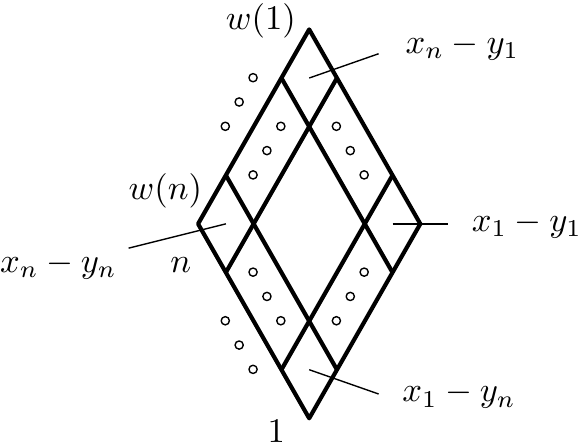}$$
\end{Th}
\noindent\textsc{Proof. }
This is clear. The puzzles are in one-to-one correspondence.
\qquad Q.E.D.
\bigbreak

Here are two examples.
$$\includegraphics{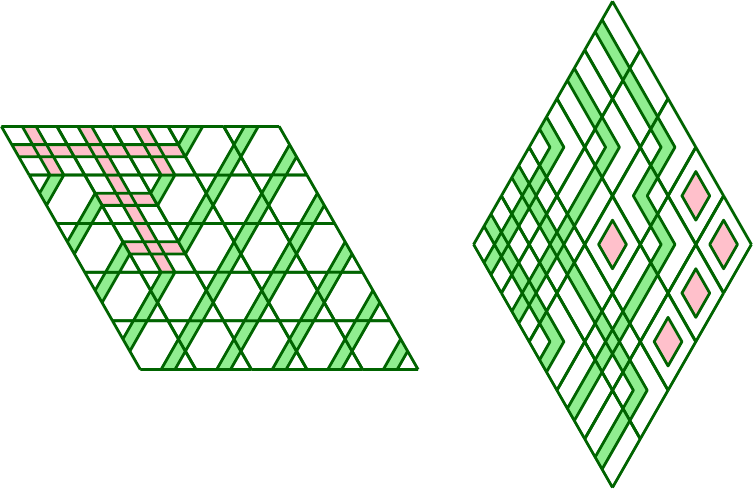}$$
We will use the mirror reflection, or rotation of the above theorem freely.

\section{Yang--Baxter Equations}

The main technique to deduce $S_w(x)=\mathfrak{S}_w$ is the so-called \emph{Yang--Baxter equations}.
The following theorem is inspired by an analogy to the Yang--Baxter equation in \cite{zinnjustin2009littlewoodrichardson}. 

\begin{Th}[Yang--Baxter Equation]\label{YBE}For any rule,
the following values of chess boards are equal
$$\includegraphics{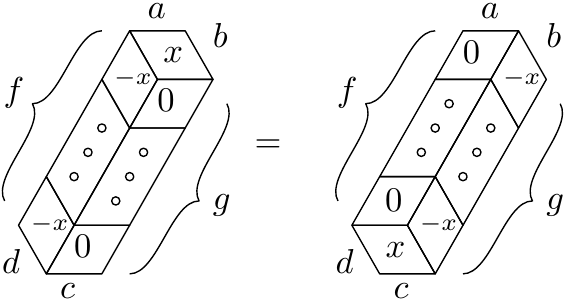}$$
where the letters around the puzzle stand for the number of pipes, where
$a+b\leq 2$ and $c+d\leq 2$; the letters inside puzzle stand for the valuation.
\end{Th}

\noindent\textsc{Proof. }
Firstly, let us prove the case for $k=1$.
This follows from case by case checking (up to rotation, and reflection)
$$\includegraphics{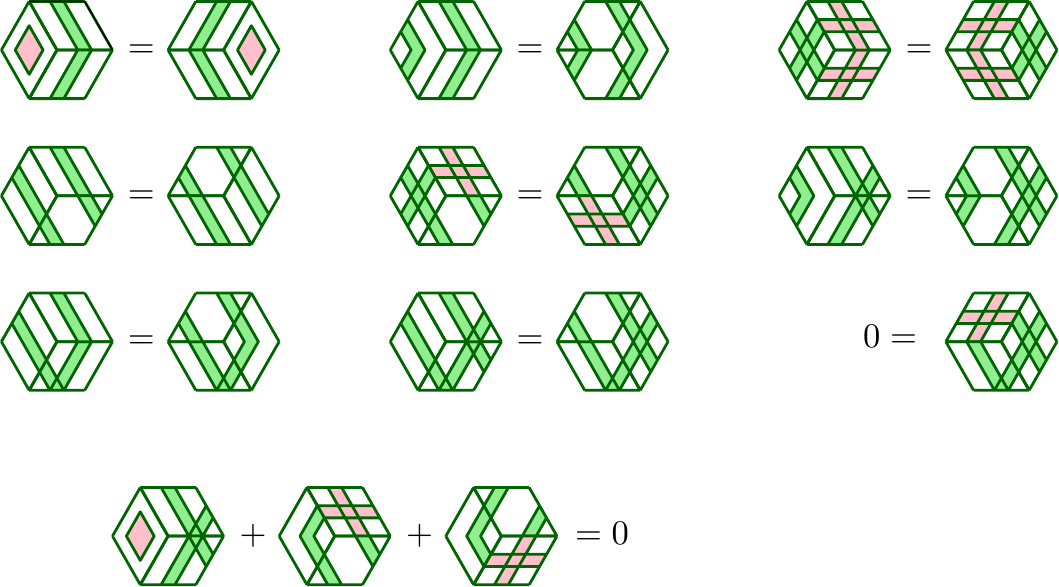}$$
Note that for the last case, we use the assumption that $x_i+y_i+z=0$.

Before processing the proof, let us denote the number of pipes around puzzle as the following
$$\includegraphics{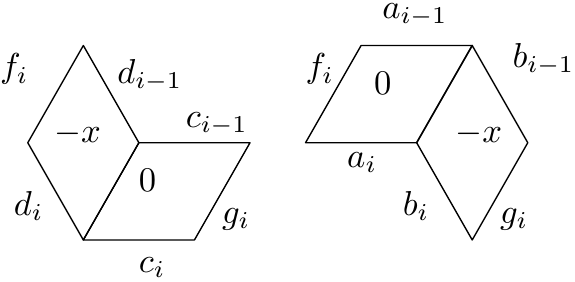}.$$
For general $k>1$, I claim that for any $I\subseteq \{1,\ldots,k-1\}$, if we denote
$$L_I=\left\{\begin{minipage}{0.3\linewidth}
solutions of the left hand side with $c_i+d_i\leq 2$ for all $i\in I$
\end{minipage}\right\},\qquad R_I=
\left\{\begin{minipage}{0.3\linewidth}
solutions of the right hand side with $a_i+b_i\leq 2$ for all $i\in I$
\end{minipage}\right\},
$$
then the sum of values of $L_I$ coincides the counterpart of $R_I$.
This follows from induction by the following diagrammatic proof
$$\includegraphics{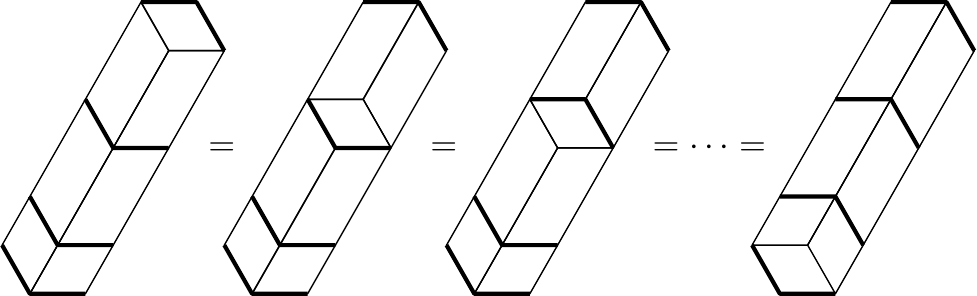},$$
where bold lines stand the restriction that the number of pipes $\leq 2$.
Note that we use the induction hypothesis and the second observation last section.
As a result,
$$L=\bigcup_{1\leq i\leq k-1} L_{\{i\}},\qquad R=\bigcup_{1\leq i\leq k-1}R_{\{i\}}$$
have the same sum of values by a simple argument of the inclusion-exclusion principle.

So it rests to deal with the case the case all $c_i+d_i=2$ for $1\leq i\leq k-1$.
The number of cases is still very limited (up to rotation)
$$\includegraphics{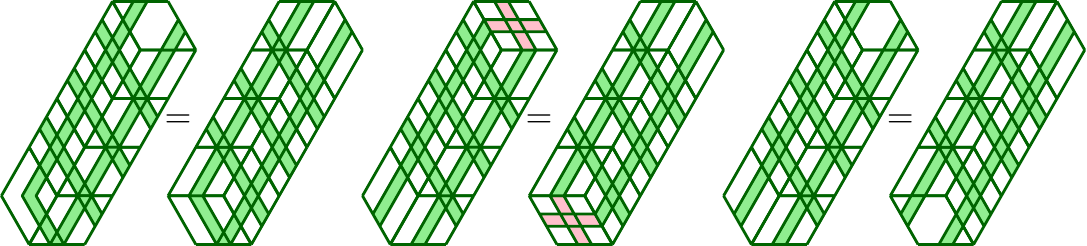}$$
This finishes the proof.
\qquad Q.E.D.
\bigbreak

Note that we cannot remove the assumption that $a+b\leq 2$ and $c+d\leq 2$. The following case is an example.
$$\includegraphics{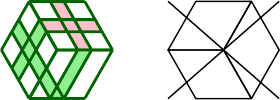}$$
Actually, this is the only counterexample (up to mirror reflection and rotation) in the hexagon of unit side length.

\section{The Puzzle Proof}

To prove Theorem \ref{ThBPD}, we need two extra lemmas to ensure in the case we considered, the conditions of Yang--Baxter equations in Theorem \ref{YBE} follows automatically.

\begin{Lemma}\label{YoungDiagram1}Assume we have the following ruled chess board.
$$\begin{array}{c|c}
\begin{minipage}{0.5\linewidth}
That is, a chess board contains 
a Young diagram at the corner as sub chess board.\par
\quad The rule is
such that single pipes $p_1,\ldots,p_m$ lie over the upper side, and no pipe on the right boundary.
\end{minipage}
&\begin{array}{c}
\includegraphics{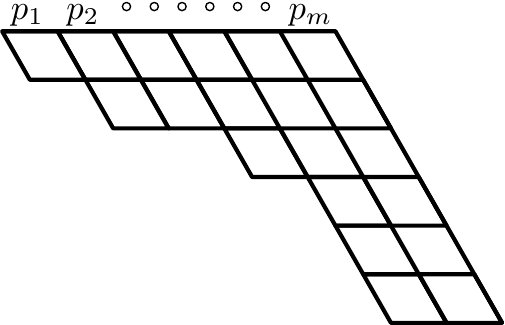}
\end{array}
\end{array}$$
Then for any solution to this rule,
only $\begin{array}{c}\includegraphics{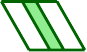}\end{array}$ will be used inside Young diaram.
\end{Lemma}
\noindent\textsc{Proof. }
The proof follows from induction from the right upper corner.
Note that for $\begin{array}{c}
\includegraphics{onehorizontalpuzzle}
\end{array}$, $(a,b)=(1,0)$ implies $(c,d)=(1,0)$ where $a,b,c,d$ are number of pipes.
\qquad Q.E.D.
\bigbreak

\begin{Lemma}\label{YoungDiagram2}Assume we have the following chess board
$$\begin{array}{c|c}
\begin{minipage}{0.5\linewidth}
That is, a chess board contains an $m\times n$ horizontal parallelogram as subshape with $m\leq n$
and a Young diagram at the corner as sub chess board.\par
\quad The rule is
such that the parallelogram is surrounded by the end points of single pipes $p_1,\ldots,p_m$ and $q_1,\ldots,q_n$.
\end{minipage}
&\begin{array}{c}
\includegraphics{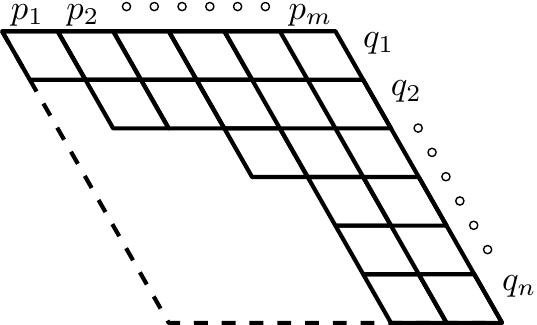}
\end{array}
\end{array}$$
Assume any pair of pipes from different families of
$$\{p_1,\ldots,p_m\},\quad \{q_1\},\quad \cdots,\quad \{q_n\}$$
does not intersect. Then for any solution to this rule,
only $\begin{array}{c}\includegraphics{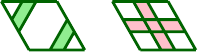}\end{array}$ will be used inside the Young diaram.
\end{Lemma}
\noindent\textsc{Proof. }
Firstly, no more than $n$ pipes could go into the parallelogram from the left wall.
Note that for any puzzle one of whose boundary parallel to the left wall has no more than $1$ pipe going through.
The only case need to check is when this wall cuts some puzzle,
$\begin{array}{c}\includegraphics{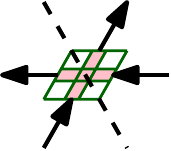}\end{array}$
but the pipe is orientated, the horizontal pipe is not into.

Secondly, since we assume $m\leq n$, so no more than $m$ pipes
could go into parallelogram from the lower wall, and they are all members of $\{q_i\}$.
Actually, by the condition of free of intersecting and $m\leq n$,
the points on pipes $q_1,\ldots,q_n$ cut by the lower wall must be in order.

As a result, there are exactly $n$ pipes going into the parallelogram from the left wall, and $m$ pipes going into from the lower wall.
To prove the conclusion,
by cutting shorter, it suffices to prove the puzzles touching the left wall is the case.
By the discussion above, for such a puzzle, the $d$ in
$\begin{array}{c}
\includegraphics{onehorizontalpuzzle}
\end{array}$ must be $1$.
Then we can prove inductively from the first row, by the fact $(a,d)=(1,1)\iff (b,c)=(1,1)$,
in which case only $\begin{array}{c}\includegraphics{RegularPuzzle}\end{array}$ will be used.
\qquad Q.E.D.
\bigbreak

\begin{Th}[=Theorem \ref{ThBPD}]For any $w\in \mathfrak{S}_n$, 
the Schubert polynomials computed by the bumpless pipe dreams $S_w(x,y)$ satisfies
$S_w(x)=\mathfrak{S}_{w}(x)$.
\end{Th}
\noindent\textsc{Proof. }
Consider the following ruled, valued chessboard.
$$\includegraphics{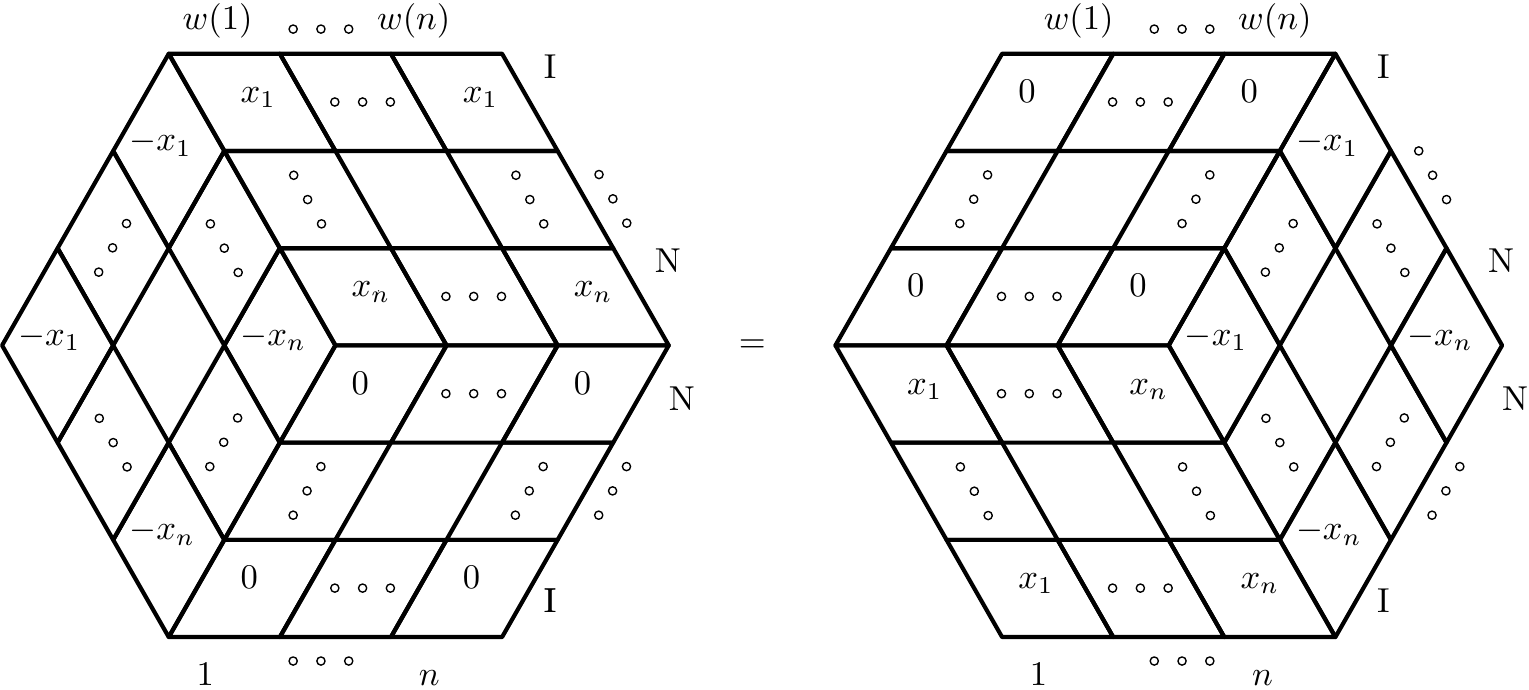}$$
I claim we can use Yang--Baxter equation in Theorem \ref{YBE} several times
to show that the value of the left-hand side equals that of the right-hand side
as the following diagram.
$$\includegraphics{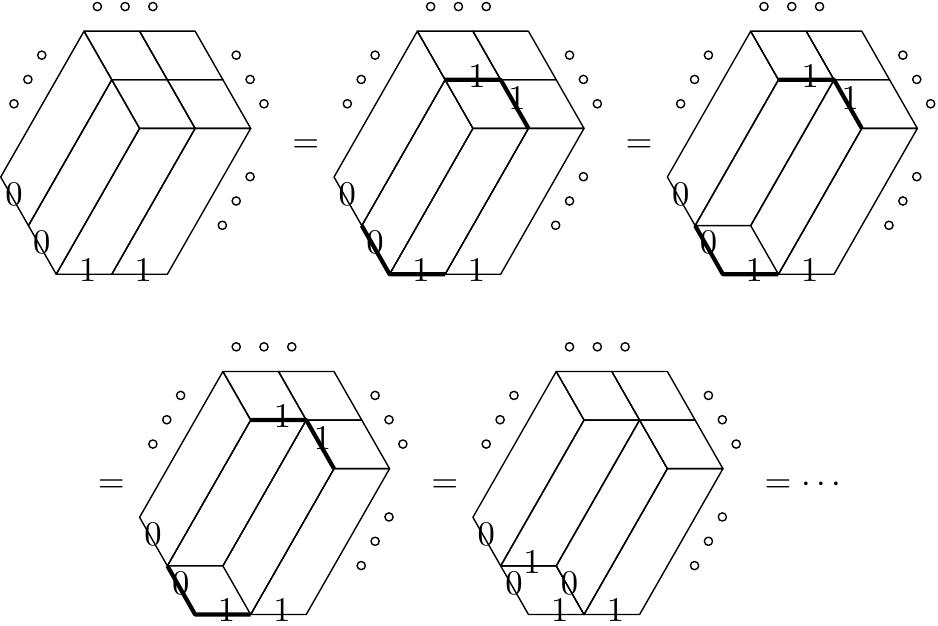}$$
Where the numbers stand the restriction on numbers of pipes.
By the two lemmas above, we can freely remove and add the restriction inside the chessboard marked inside the diagram.

The value of left hand side is
$$\sum_{w\stackrel{\text{reduced}}=u\cdot v}(-1)^{\ell(v)}S_v(x)\mathfrak{S}_u(x)$$
due to Thereom \ref{puzzlePD}, Theorem \ref{puzzleBPD} and Lemma \ref{YoungDiagram2}.
The value of right hand side is $\begin{cases}1,&w=\operatorname{id},\\
0,&\text{otherwise},
\end{cases}$ by Lemma \ref{YoungDiagram1}.
But by the Lemma below this characterizes $\mathfrak{S}_w(x)$.
\qquad Q.E.D.
\bigbreak

\begin{Lemma}\label{ChofSch}Let $\{S_w(x)\}$ be a series of polynomials parameterized by permutation $w\in \mathfrak{S}_\infty$.
If for any $w\in \mathfrak{S}_\infty$,
$$\sum_{w\stackrel{\text{reduced}}=u\cdot v}(-1)^{\ell(v)}S_{v^{-1}}(x)\mathfrak{S}_{u}(x)=\begin{cases}
1,& w=\operatorname{id}, \\
0,& \text{otherwise}.
\end{cases}$$
then $S_w(x)=\mathfrak{S}_w(x)$.
\end{Lemma}
\noindent\textsc{Proof. }%
Note that when $S_w(x)=\mathfrak{S}_w(x)$, the equality holds, see \cite{lam2018stable}. 
Conversely, it suffices to show that for a series of polynomials $\{\varphi_w(x)\}$, if for all $w\in \mathfrak{S}_\infty$
$$\sum_{w\stackrel{\text{reduced}}=u\cdot v}(-1)^{\ell(v)}\varphi_{v^{-1}}(x)\mathfrak{S}_{u}(x)=0,$$
then $\varphi_w=0$.
This follows from the fact $\mathfrak{S}_{\operatorname{id}}(x)=1$ and an easy argument of induction.
\qquad Q.E.D.
\bigbreak

This finishes the proof. 

\section{Remarks}

In this section, I will give some remarks on the difficulty of applying this method on
double Schubert polynomials. 
Please contact me (see last page) without any hesitation if the reader figures out a proof for the Double Schubert polynomials.

One may think to modify Yang-Baxter equation in Theorem \ref{YBE} by
$$\includegraphics{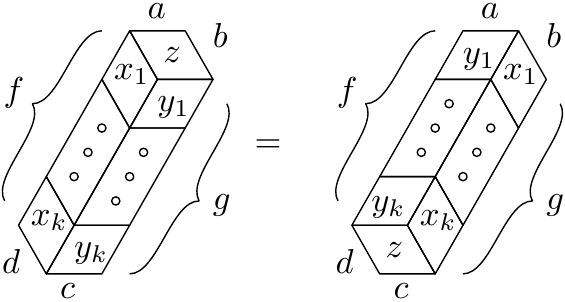}$$
whenever $x_i+y_i=z$ for all $i$ and consider the chessboard like this
$$\includegraphics{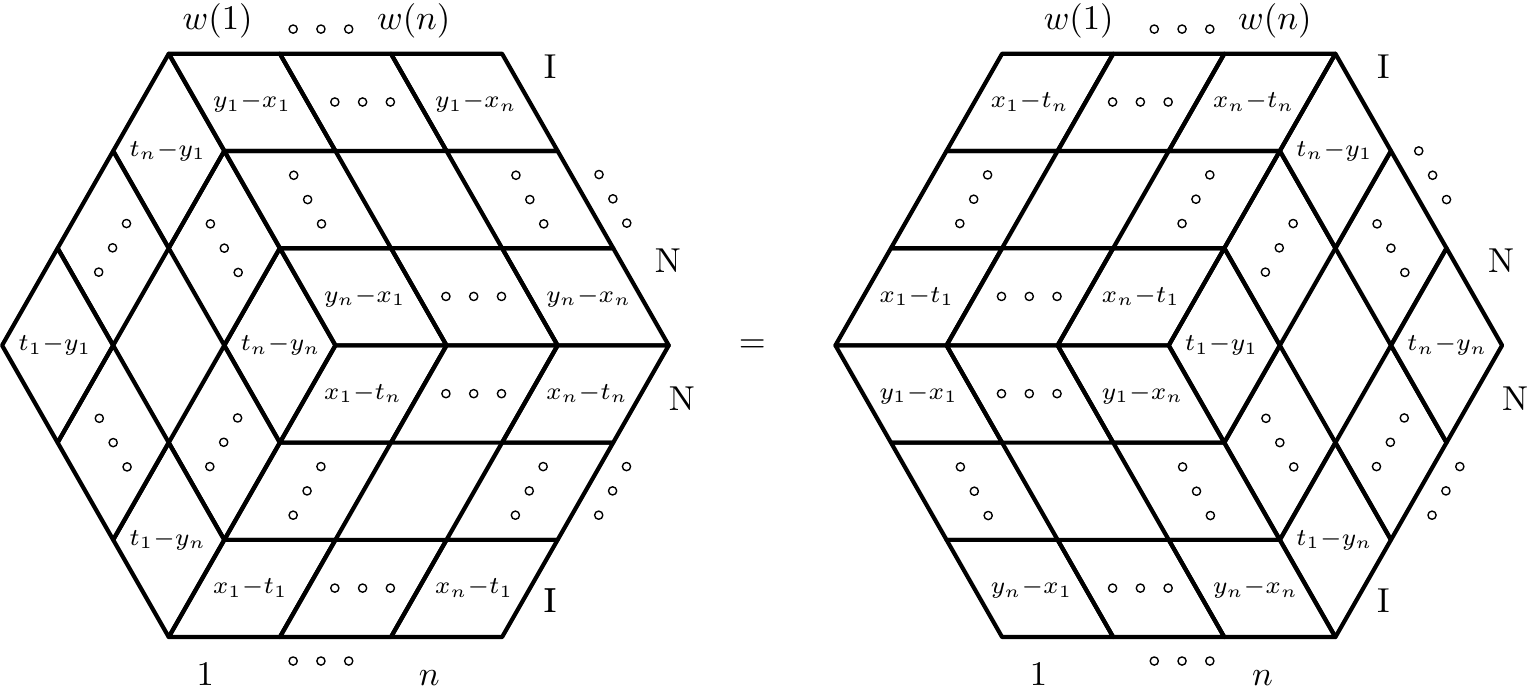}$$

The Yang-Baxter equation does not hold. 
Here are the only two counterexamples in the hexagon of unit side length.
$$\includegraphics{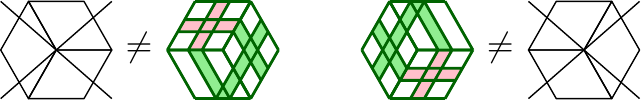}. $$
Note that the equalities hold only for $y_i=0$, which is ensured in Theorem \ref{YBE}.

Nevertheless, the values of two chessboards literally equal. 
Since it is proved in \cite{lam2018stable} that $S(x,t)=\mathfrak{S}(x,t)$, 
the value of left hand side is 
$$\sum_{w=a\cdot b\cdot c}\mathfrak{S}_{c}(x,t)S_b(t,y)\mathfrak{S}_a(y,x)=
\sum_{w=a\cdot b\cdot c}\mathfrak{S}_{c}(x,t)\mathfrak{S}_b(t,y)\mathfrak{S}_a(y,x)=\mathfrak{S}_{w}(x,x)=\begin{cases}1,&w=\operatorname{id},\\
0,&\text{otherwise},
\end{cases}$$
exactly the value of the right hand side. 
Of course, if one can prove the equality of values, then this will characterize the double Schubert polynomials. 

But the author has not figured out how to prove it. 
The following example indicates that that cleverly cancel out.
$$\includegraphics{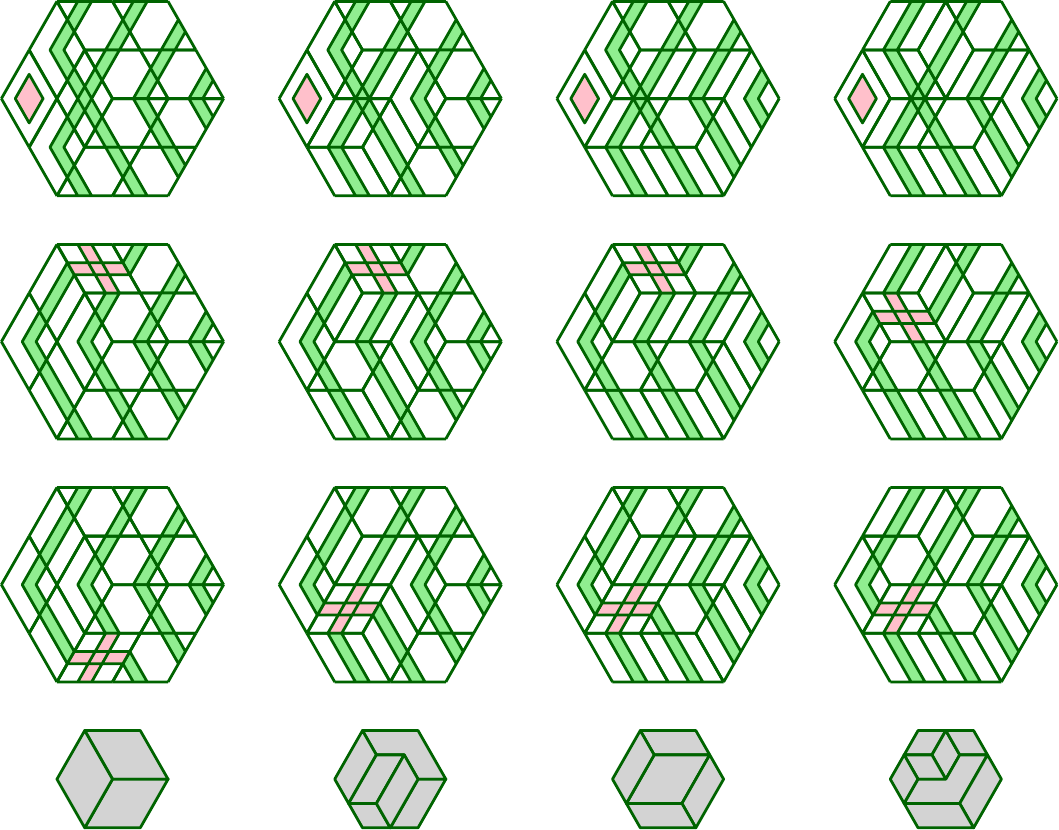}$$

Another attempt is to show that $S_{w}(x,w'x)=0$ for all $\ell(w')<\ell(w)$. 
This is sufficient to force $S_w(x,y)=\mathfrak{S}_w(x,y)$. 
Actually, we have 

\begin{Lemma}\label{ChofSch2}For a series of polynomials $\{S_w(x,y)\}$ parameterized by permutation $w\in \mathfrak{S}_\infty$,
if $\deg S_w=\ell(w)$, $S_w(x,0)= \mathfrak{S}_w(x)$ and
$$\ell(w')<\ell(w)\Longrightarrow S_w(x,w'x)= 0, $$
then $\mathfrak{S}_w(x,y)=S_w(x,y)$.
\end{Lemma}
\noindent\textsc{Proof. }%
We can expand $S_w(x,y)$ in the form $\sum_{u\in \mathfrak{S}_\infty} c^u(y)\mathfrak{S}_u(x,y)$.
Then taking in $y=0$ and by degree reason, it suffices to show $c^{u}(y)=0$ for all $u$ with $\ell(u)<\ell(w)$.
Otherwise, pick the minimal $u$ such that $c^{u}(y)\neq 0$,
then $0=S_{w}(x,ux)=c^{u}(ux)\mathfrak{S}_u(x,ux)$, implies $c^u(ux)=0$, i.e. $c^u=0$.
\qquad Q.E.D.
\bigbreak

From the bumpless pipe dream, it is easy to show that $S_w(x,w'x)=0$ when $w'< w$ in Bruhat order. 
But in general, the author do not figure out a proof. 
%

\bibliographystyle{plain}
\bibliography{bibfile}

\begin{thebibliography}{Knu19}

\bibitem[BB93]{bergeron1993rc}
Nantel Bergeron and Sara Billey.
\newblock Rc-graphs and schubert polynomials.
\newblock {\em Experimental Mathematics}, 2(4):257--269, 1993.

\bibitem[FK96]{fomin1996yang}
Sergey Fomin and Anatol~N Kirillov.
\newblock The yang-baxter equation, symmetric functions, and schubert
  polynomials.
\newblock {\em Discrete Mathematics}, 153(1-3):123--143, 1996.

\bibitem[FS94]{fomin1994schubert}
Sergey Fomin and Richard~P Stanley.
\newblock Schubert polynomials and the nilcoxeter algebra.
\newblock {\em Advances in Mathematics}, 103(2):196--207, 1994.

\bibitem[Knu19]{knutson2019schubert}
Allen Knutson.
\newblock Schubert polynomials, pipe dreams, equivariant classes, and a
  co-transition formula, 2019.

\bibitem[LLS18]{lam2018stable}
Thomas Lam, Seung~Jin Lee, and Mark Shimozono.
\newblock Back stable schubert calculus, 2018.

\bibitem[LS82]{lascoux1982structure}
Alain Lascoux and Marcel-Paul Sch{\"u}tzenberger.
\newblock Structure de hopf de l'anneau de cohomologie et de l'anneau de
  grothendieck d'une vari{\'e}t{\'e} de drapeaux.
\newblock {\em CR Acad. Sci. Paris S{\'e}r. I Math}, 295(11):629--633, 1982.

\bibitem[ZJ09]{zinnjustin2009littlewoodrichardson}
P.~Zinn-Justin.
\newblock Littlewood--richardson coefficients and integrable tilings, 2009.

\end{thebibliography}
\noindent
\hfill\rule{0.6\linewidth}{0.4pt}\hfill\mbox{}
\begin{center}
XIONG Rui, master stundent\smallbreak
\textit{\href{http://english.spbu.ru/}{Saint Petersburg State University}}
\textit{\href{https://math-cs.spbu.ru/en/}{Department of Mathematics and Computer Science}}
Saint Petersburg, 199178, Russia, Line 14th (Vasilyevsky Island),
\par
e-mail: \url{XiongRui_Math@126.com},\par
homepage: \url{www.cnblogs.com/XiongRuiMath}.
\end{center}

\end{document}